\documentclass[12pt,a4paper]{article}
\usepackage{amsmath, amssymb, theorem, latexsym}
\usepackage{graphicx}

\newtheorem{theorem}{Theorem}[section]
\newtheorem{lemma}[theorem]{Lemma}
\newtheorem{proposition}[theorem]{Proposition}
\newtheorem{definition}[theorem]{Definition}

\newtheorem{conjecture}[theorem]{Conjecture}
\newtheorem{corollary}[theorem]{Corollary}

\newtheorem{remark}[theorem]{Remark}

\allowdisplaybreaks

\def\neweq#1{\begin{equation}\label{#1}}
\def\endeq{\end{equation}}
\newcommand\finbox{~\hfill$\Box$}%

\def\Om {{\Omega}}
\def\la {{\lambda}}

\def \Inte{{\rm Int\,}}

\newcommand{\beq}{\begin{equation}}
\newcommand{\eeq}{\end{equation}}
\newcommand {\ar}{\rightarrow}
\newcommand {\pa}{\partial}

\numberwithin{equation}{section}
\begin{document}
{\centering
\bfseries
{\Large On spectral  minimal partitions~:  the case of the sphere}

\par
\mdseries
\scshape
\small
B. Helffer$^*$\\
T. Hoffmann-Ostenhof$^{**,***}$ \\
S. Terracini$^{****}$.
\par
\upshape
D\a'epartement de Math\a'ematiques, Univ Paris-Sud and CNRS $^*$\\
Institut f\"ur Theoretische Chemie, Universit\"at Wien$^{**}$\\
International Erwin Schr\"odinger Institute for Mathematical
Physics$^{***}$\\
Universita di Milano Bicocca$^{****}$.\\
\today

}

\begin{abstract}
In continuation of \cite{HHOT}, \cite{HHO}, \cite{BHV:2007} 
 and \cite{BHHO}, we analyze   the
properties of  spectral minimal  partitions and focus  in this
paper our
analysis on  the case of  the sphere. We  
 prove that a minimal $3$-partition for the sphere
 $\mathbb S^2$ 
 is  up to rotation the so called  ${\bf Y}$-partition. This
 question
 is connected to  a celebrated conjecture of Bishop in harmonic analysis. 
\end{abstract}


\section{Introduction}
Motivated by questions related to some conjecture of Bishop \cite{Bis}, we
continue
 the analysis of spectral minimal partitions developed for plane
 domains
 in  \cite{HHOT}, \cite{HHO}, \cite{BHV:2007}, \cite{BHHO} and analyze
 the case of the two-dimensional sphere $\mathbb S^2$.\\

In the whole  paper the Laplacian is the Laplace-Beltrami operator on
$\mathbb S^2$. 
We describe as usual $\mathbb S^2$ in $\mathbb R^3_{x,y,z}$ by the spherical  coordinates,
 \begin{equation}
x= \cos \phi\sin \theta , y =\sin \phi \sin \theta, z= \cos
\theta\,,\,
 \mbox{ with } \phi\in [-\pi,\pi[\,, \, \theta \in ]0,\pi[\,,
\end{equation}
and we add the two
 poles ``North'' and ``South'', corresponding to the two points 
  $(0,0,1)$ and $(0,0,-1)$. \\
If $\Omega$ is a regular bounded open set 
 with piecewise $C^{1,+}$ boundary\footnote{ i.e. with piecewise 
$C^{1,\alpha}$ boundary, for some $\alpha >0$},
 we consider the Dirichlet Laplacian $H=H(\Omega)$ and we would like to analyze
 the question of the existence and of the properties of minimal
 partitions
 for open sets $\Omega$ on $\mathbb S^2$. When $\Omega$ is strictly
 contained
 in $\mathbb S^2$, the question is not fundamentally 
different of the case of planar domains, hence 
 we will focus on the case of the whole sphere $\mathbb S^2$ and 
 on the search for  possible candidates for minimal $k$-partitions 
 of the sphere for $k$ small.\\

To be more precise, let us now recall a few definitions that the reader
can for example find in \cite{HHOT}. 
For  $1\le k\in \mathbb N$ and $\Omega\subset \mathbb S^2$, we call  a
{\bf spectral}
$ k$-\textbf{partition}\footnote{We say more shortly ``$k$-partition''.} of $\Om$   a family $\mathcal D=\{D_i\}_{i=1}^k$ of pairwise  disjoint 
open regular domains  such that 
\begin{equation}\label{opar}
\cup_{i=1}^kD_i\subset \Om. 
\end{equation}
It is called  \textbf{strong} if 
\begin{equation}\label{regcov}
\Inte(\overline{\cup_{i=1}^k D_i})\setminus \pa \Om=\Om\,.
\end{equation}
We denote by $\mathfrak O_k$ the set of such partitions.
 For $\mathcal D\in \mathfrak O_k$ we introduce 
\begin{equation}\label{LA}
\Lambda(\mathcal D)=\max_{i}\la(D_i)\,,
\end{equation}
where $\lambda(D_i)$ is the ground state energy of $H(D_i)\,$, 
and 
\begin{equation}\label{Lfrak}
\mathfrak L_k(\Omega) =\inf_{\mathcal D\in \mathfrak O_k}\Lambda(\mathcal D)\,.
\end{equation}
We call  a spectral minimal  $k$-partition, a $k$-partition 
$\mathcal D\in \mathfrak O_k$ 
such that  $$  \mathfrak L_k(\Omega)=\Lambda(\mathcal
 D)\,.
$$
More generally we can consider (see in \cite{HHOT}) for $p\in [1,+\infty[$
\begin{equation}\label{LAp}
\Lambda^p(\mathcal D)=(\frac 1k \sum_i \la(D_i)^p)^{\frac 1p}\,,
\end{equation}
and
\begin{equation}\label{Lfrakp}
\mathfrak L_{k,p}(\Omega)=\inf_{\mathcal D\in \mathfrak O_k}\Lambda^p(\mathcal D)\,.
\end{equation}
We write  $\mathfrak L_{k,\infty}(\Omega)=\mathfrak L_k(\Omega)\,$ and
recall the monotonicity property
\beq\label{monoto}
\mathfrak L_{k,p} (\Omega) \leq \mathfrak L_{k,q}(\Omega)\,\mbox{ if }
p\leq q\,.
\eeq 
The notion of $p$-minimal $k$-partition can be extended accordingly,
by minimizing $\Lambda^p(\mathcal D)$.\\
We would like to give in this article the proof of the following theorem.
\begin{theorem}\label{conj1}~\\
Any minimal  $3$-partition of $\mathbb S^2$ is up to a fixed rotation
obtained by the so called {\bf Y}-partition whose boundary is given by the intersection
of $\mathbb S^2$ with  the three half-planes  defined respectively by
$\phi =0, \frac{2\pi}{3}, \frac{-2 \pi}{3}\,$. Hence 
\begin{equation}\label{valeur}
\mathfrak L_3(\mathbb S^2)=\frac{15}{4}\,.
\end{equation}
\end{theorem}
This theorem is immediately related (actually a consequence of) to
a conjecture of Bishop  (Conjecture 6)  proposed in \cite{Bis} stating that~:
\begin{conjecture}[Bishop 1992]\label{conj2}~\\
The minimal $3$-partition for $\frac 13 (\sum_{i=1}^3 \lambda(D_i))$
 corresponds to the {\bf Y}-partition.
\end{conjecture}
We can indeed observe that if for some $(k,p)$ there exists a  $p$-minimal
$k$-partition $\mathcal D_{k,p}$ such that $\lambda(D_i) =\lambda(D_j)$ for all $i,j$,
then by the monotonicity property $\mathcal D_{k,p}$ is a $q$-minimal
partition for any $q\geq p$.\\
\begin{remark}~\\
At the origin, Bishop's Conjecture was motivated by the analysis
 of the properties of Harmonic functions in conic sets. The whole
 paper by Friedland-Hayman \cite{FrHa} (see also references therein) which inspires our Section \ref{s5} is written
 in this context. The link between our problem of minimal partitions
 and the problem in harmonic analysis can be summarized in this way.
If we consider a homogeneous Lipschitzian function of the form 
 $u(x)=r^\alpha g(\theta,\phi)$ in $\mathbb R^3$,
 which is harmonic outside its nodal set
 and such that the complementary of the nodal set divides
 the sphere in three parts, then
 $$ \alpha (\alpha +1)\geq \mathfrak L_3(\mathbb S^{2}) \,.$$
Hence Theorem \ref{conj1} (and more
specifically
 \eqref{valeur}) implies  $\alpha\geq 3/2 $.
This kind of property can be
 useful to improve some statements in \cite{CL1,CL2} (see inside the
 proofs of Lemmas~2 in
 \cite{CL1} and 4.1 in \cite{CL2}). 
\end{remark}
A similar question was analyzed (with partial success) when looking in
\cite{HHO} at
candidates of minimal $3$-partitions  of the unit disk $D(0,1)$ in
$\mathbb R^2$. The most natural
candidate
 was indeed  the Mercedes Star, which is  the $3$-partition 
given   by three disjoint sectors with opening 
angle $2\pi/3$, i.e. 
 \begin{equation}\label{sector}
D_1=\{x\in \Om\:|\:\omega\in ]0,2\pi/3[\}
\end{equation}
and $D_2, D_3$ are obtained by rotating $D_1$ by $2\pi/3\,$,
respectively by  $4\pi/3\,$.
Hence the Mercedes star in \cite{HHO}  is replaced here by the ${\bf
  Y}$-partition in Theorem \ref{conj1}. 
We observe that  {\bf Y}-partition  can also be described  the inverse image of the mercedes-star
 partition  by the map
$\mathbb S^2\ni (x,y,z)\mapsto (x,y)\in D(0,1)$.\\

Here let us mention the two main statements giving the proof of Theorem~\ref{conj1}.
\begin{proposition}\label{antipodal}~\\
If  $\mathcal D=(D_1,D_2,D_3)$ is a $3$-minimal partition, then its
boundary
 contains two antipodal points.
\end{proposition}
The proof of Proposition \ref{antipodal} will be achieved in Section \ref{s6} and involves Euler's
formula
 and the thorem of  Lyusternik and Shnirelman.
\begin{proposition}\label{centera}~\\
If there exists a  minimal $3$-partition $\mathcal D=(D_1,D_2,D_3)$ 
 of $\, \mathbb S^2$ with two antipodal points in 
$\cup_{i} \pa D_i$,  then it is (after possibly a rotation) 
 the  ${\bf Y}$-partition.
\end{proposition}
The proof of Proposition \ref{centera}  will be done in  Section \ref{s8} by lifting this  $3$-partition 
  on the double covering $\mathbb S^2_{\mathcal C}$ of $\ddot{\mathbb S^2}$, where
  $\ddot{\mathbb S^2}$ is the
  sphere 
 minus two antipodal points.\\

More precisely, following what has been done in the approach of the Mercedes-star
conjecture in \cite{HHO}, the steps for the proof 
 of Theorem \ref{conj1} (or towards Conjecture \ref{conj2} if we were
 able to 
 show that the minimal partition for $\Lambda^1$ has all the $\lambda(D_j)$  equal)  
are the following~:
\begin{enumerate}
\item
One has to prove that minimal partitions on $\mathbb S^2$ exist
 and share the same properties as for planar domains : regularity and 
 equal  angle meeting  property. This will be done in Section~\ref{s2}.
\item One can observe that  the
 minimal $3$-partition cannot be a nodal partition. This is a
 consequence of Theorem \ref{L=L} in Section \ref{s2} and  
  of the fact that the multiplicity of the second
  eigenvalue (i.e. the first non zero one) is more than $2$ actually
   $3$. 
\item The Euler formula implies that there exists only one
 possible  type of minimal $3$-partitions. Its boundary  
 consists of  two points $x_1$ and $x_2$ and three arcs
 joining  these two points. This will be deduced in Subsection~\ref{s3}.
\item The next point is  to show a  minimal partition
  has in its boundary two antipodal points. 
\item The next point is that any minimal $3$-partition  which contains
  two antipodal  points in its boundary can be lifted in a symmetric $6$-partition
  on the double covering $\mathbb S^2_{\mathcal C}$. More precisely, 
if $\mathcal D
 =(D_1,D_2,D_3)$
 and $\Pi$ is the canonical projection of $\mathbb S^2_{\mathcal C}$
 onto $\ddot{\mathbb S^2}$ , we get the $6$-partition $\mathcal
   D_{\mathcal C}$
 by considering 
$$
 \mathcal
   D_{\mathcal C} = (D_1^+,D_2^+,D_3^+, D_1^{-}, D_2^{-}, D_3^-)\,,
$$
where for $j=1,2,3$ $D_j^{+}$ and $D_{j}^{-}$ denote the two
components
 of $\Pi^{-1} (D_j)$. If $\mathcal I$ denotes the map on $\mathbb
 S^2_{\mathcal C}$
 defined, for $m\in \mathbb S^2_{\mathcal C}$  by
\begin{equation}\label{defIrev}
\Pi (\mathcal I(m))= \Pi(m) \mbox{ with }\mathcal I(m)\neq m \,,
\end{equation}
we observe that
$$
\mathcal I (D_j^+)=D_j^-\;.
$$

\item The last point is to show that on this  double covering
 a minimal  symmetric $6$-partition is necessarily the double {\bf
 Y}-partition, which is the inverse image in $\mathbb
 S^2_{\mathcal C}$
 of the {\bf Y}-partition.   
\end{enumerate}
All these points will be detailed in the following sections
 together with analogous questions in the case of minimal
 $4$-partitions.\\

In the last section, we will describe what can be said towards
 the proof of Bishop's conjecture and about the large $k$ behavior
 of $\mathfrak L_k$ using mainly the tricky estimates of
 Friedland-Hayman
 \cite{FrHa}.

\section{Definitions, notations and extension of  previous results to
  the sphere.}\label{s2}
We first recall more  notation, definitions and results essentially 
  extracted
 of \cite{HHOT} but which have to be extended from the case of planar
 domains  to the case of domains  in $\mathbb S^2$. 
\begin{definition}\label{def:regopenset}~\\
We say that an open  domain  $D\subset  \mathbb S^2$ is
 \textbf{regular}
 if it satisfies an interior cone condition and if 
$\pa D$ is the union a finite number of simple regular arcs 
$\gamma_i(\overline{I_i})$, with $\gamma_i\in\mathcal C^{1,+}(\overline I_i)$ with no mutual 
nor self intersections, except possibly at the endpoints.
\end{definition}

For a given set $\Omega\subset \mathbb S^2$, we  are
interested in the eigenvalue problem for $H(\Om)$, the Dirichlet
realization of the Laplace Beltrami operator in $\Omega$.  We shall denote for any open domain
$\Omega$
 by $\lambda(\Omega)$  the lowest 
eigenvalue of $H(\Omega)$. 
 We define for any eigenfunction $u$ of $H(\Omega)$ 
\begin{equation}
N(u)=\overline{\{x\in \Om\:\big|\: u(x)=0\}}
\end{equation}
and call the components of $\Om\setminus N(u)$ the nodal domains of $u$. The number of 
nodal domains of such a function will be called $\mu(u)$. 

If $\mathcal D$ is a strong partition, 
we say $D_i,D_j\in \mathcal D$ are \textbf{neighbors}
if 
\begin{equation}\label{DsimD}
\Inte(\overline{D_i\cup D_j})\setminus \pa \Om\text{ is connected}
\end{equation}
and write in this case  $D_i\sim D_j$.  
We  then define   a {\bf graph}  $G(\mathcal D)$
by associating to
 each $D_i\in \mathcal D$  a vertex $v_i$ and to  each pair $D_i\sim D_j$ we associate 
an edge $e_{i,j}$. 

Attached to a regular partition $\mathcal D$ we can associate its
boundary  $N=N(\mathcal D)$ which is the closed set in $\overline\Om$
defined by 
\begin{equation}
N(\mathcal D)=\overline{\bigcup_i(\pa D_i\cap \Om)}.
\end{equation}
This leads us to introduce the  set $\mathcal M(\Om)$ of the regular closed sets.

\begin{definition}\label{AMS}~\\
A closed set $N\subset \overline\Om$ belongs to $\mathcal M(\Om)$ if $N$ meets the following requirements:
\\
\textbf{(i)} 
There are finitely many distinct critical points $x_i\in \Om\cap N$ and associated positive 
integers $\nu(x_i)$  with $\nu(x_i)\ge 3$
such that, in a sufficiently small neighborhood of each of the $x_i$, 
$N$ is the union of $\nu(x_i)$ disjoint 
 (away from $x_i$ non self-crossing) smooth arcs  with one end at $x_i$ (and each pair 
defining at $x_i$ a positive angle in $]0,2\pi[$)
and such that in the complement of these points in $\Om$, $N$ 
is locally diffeomorphic to a smooth arc. We denote by $X(N)$ the set of these critical points.  \\
\textbf{(ii)} 
$\pa\Om\cap N$ consists of a (possibly empty) finite set of points $z_i$, such that at each $z_i$, 
$\rho(z_i), \text{with } \rho(z_i)\ge 1$  arcs hit the boundary. Moreover for each $z_i\in \pa \Om$,  
then $N$ is near $z_i$ the union of $\rho_i$ distinct smooth arcs  which hit 
$z_i$ with strictly positive distinct angles.  We denote by $Y(N)$ the set of these critical points. \end{definition}

Conversely, if $N$ is a regular closed set, then the family $\mathcal D(N)$ of connected 
components of $\Om\setminus N$ belongs (by definition) to $\mathcal R(\Om)$, hence 
regular and strong.

\begin{definition}\label{eacp}~\\
We will say that a closed set has the \textbf{equal angle meeting
  property} ({\bf eamp}), if the arcs meet 
with equal angles at each critical point  $x_i\in N\cap \Om$ and also with equal angles 
at the $z_i\in N\cap  \pa \Om$. For the boundary points $z_i$ we mean that the two arcs in the 
boundary are included.
\end{definition}

We will say that the partition is {\bf eamp}-regular if it is regular and
satisfies the equal angle meeting property.\\

It has been proved by Conti-Terracini-Verzini \cite{CTV0, CTV2, CTV:2005} that 
\begin{theorem}\label{thstrreg}~\\
For any $k$ there exists a minimal {\bf eamp}-regular strong  $k$-partition.
\end{theorem}
It has been proved in \cite{HHOT} the 
\begin{theorem}\label{thanyreg}~\\
Any  minimal spectral $k$-partition admits a representative which is 
{\bf eamp}-regular and strong.
\end{theorem}

A basic result concerns  the regularity (up to the boundary if any) of the nodal partition associated to an eigenfunction.

We first observe that the  results about minimal partitions for 
 plane domains can be transfered to the sphere $\mathbb S^2$. 
It is indeed enough
 to use the stereographic projection on the plane which gives  an
 elliptic operator on the plane with analytic coefficients.This map is
 a conformal map, hence respecting the angles. The
 regularity questions being local
 there are no particular problem for
 recovering the equal angle meeting  property.

A natural question is whether a minimal partition is the nodal partition 
induced by an eigenfunction.
The next  theorem gives a simple criterion for a partition to be 
 associated to a nodal set.  For this we need some additional
 definitions.

 We  recall  that the graph $G(\mathcal D)$ is \textbf{bipartite} if its vertices 
can be colored by two colors (two neighbours having  different
colors). 
In this case, we say that the partition is \textbf{admissible}.  We recall that  a collection of nodal domains of an eigenfunction is always
 admissible.

We have now the following converse theorem \cite{HHOT}~:
\begin{theorem}\label{partnod}~\\
An admissible  minimal  $k$-partition is nodal, i.e. associated to the nodal set of an 
eigenfunction of $H(\Omega)$ corresponding to an eigenvalue
 equal to $\mathfrak L_k(\Omega)$.
\end{theorem}

This theorem was already obtained for planar domains in \cite{HH:2005a} by adding a strong
a priori regularity 
and the assumption that $\Omega$ is simply
connected. Any subpartition of cardinality $2$ corresponds indeed
 to a second eigenvalue and the criterion of pair compatibility (see
 \cite{HH:2005a}) can be applied.\\

A natural question is now to determine how general is  the situation
described in Theorem \ref{partnod}. As for partitions in planar
domains,
 this can only occur in very particular cases when $k>2$.
If $\lambda_k(\Omega)$ denotes the $k$-th eigenvalue of the Dirichlet
realization of the Laplacian in an open set $\Omega$ of $\mathbb S^2$, 
the Courant Theorem says~:
\begin{theorem}\label{ThC}~\\
 The number of nodal domains $\mu(u)$ of
 an eigenfunction $u$ associated with $\lambda_k(\Omega)$ 
satisfies
 $\mu(u)\leq k\,.
$
\end{theorem}
 Then we say, as in 
\cite{HHOT}, that  $u$ is \textbf{Courant-sharp} if $\mu(u)=k$.
For any integer $k\ge 1$, we denote by $L_k(\Omega)$ the smallest eigenvalue whose eigenspace 
contains an eigenfunction with $k$ nodal domains. In general we have
\begin{equation}
\lambda_k(\Omega)\leq\mathfrak L_k(\Omega) \leq L_k(\Omega) \,.
\end{equation}

The next result of \cite{HHOT} gives  the full picture of the equality cases~:

\begin{theorem}\label{L=L}~\\
Suppose that $\Om\subset \mathbb S^2$ is regular. 
If $\mathfrak L_k(\Omega) =L_k(\Omega)$ or $\lambda_k(\Omega) =\mathfrak L_k(\Omega)$, then $$\la_k(\Omega)=\mathfrak L_k(\Omega)=L_k(\Omega)\,,$$
and  any minimal $k$-partition is nodal and admits a representative
  which is the family of nodal domains of some eigenfunction $u$
associated  to $\la_k(\Omega)$. 
\end{theorem}

This theorem will be quite useful for showing for example that, for $k=3$ and
$k=4$, a $k$-minimal partition of $\mathbb S^2$ cannot be nodal. This
will be further discussed in Section \ref{s4} (see Theorem \ref{Unter}).

\section{Courant's nodal Theorem with inversion symmetry.}\label{s4}
We collect here some easy useful observations  for the
analysis of the sphere. These considerations already appear in
\cite{Ley1} but the application to minimal partitions is new.  
We consider Courant's Nodal Theorem for $
H(\Om)$
where $\Om$ is an open connected set in $\mathbb S^2$. 
Let 
\begin{equation}
\mathbb S^2\ni (x,y,z)\mapsto I(x,y,z)=(-x,- y,-z)
\end{equation}
 denotes the  inversion map
  and 
assume that  
\begin{equation}\label{OVinv}   
I\Om =\Om \,.
\end{equation}
Note that $\Omega =\mathbb S^2$ satisfies the condition.
These assumptions imply that we can write $H(\Om)$
as a direct sum 
\begin{equation}\label{HSHA} 
H(\Om)=H_S(\Om)\bigoplus H_A(\Om)\,,
\end{equation}
where $H_S(\Om)$ and $ H_A(\Om)$ are respectively
the restrictions of $ H(\Om)$ to the $I$-symmetric 
(resp. antisymmetric)
 $L^2$-functions in $\Omega$ in $D( H(\Om))$.

For simplicity we just write $H_S,H_A$.
For the spectrum of $H(\Om)$, $\sigma$, we have 
\begin{equation*}
\sigma=\sigma_S\cup \sigma_A
\end{equation*}
so that  $\sigma_S=\{\la^S_k\}_{k=1}^\infty$ and analogously $\sigma_A=\{\la_k^A\}_{k=1}^\infty$. It is of independent 
interest to investigate how $\sigma_S$ and $\sigma_A$ are related. Obviously we have 
$$
\la_1^S<\la_1^A\le \la_2^A\mbox{ and } \la_1^S<\la_2^S\,,
$$ by standard spectral  theory.

In the present situation we can ask the question 
of a theorem \`a la Courant separately for the eigenfunctions 
of $H_S$ and $H_A$.

First we  note the following easy properties~:
\begin{enumerate}
\item 
Suppose $u$ is an eigenfunction and that $u$ is either symmetric or
anti-symmetric. Then $I\, N(u)=N(u)$, i.e. the 
nodal set is symmetric with respect to inversion.
\item 
If  $u^A$ is an eigenfunction of $H_A$, then for each nodal domain
$D_i$ of $u^A$, $ID_i$ is a distinct  nodal domain of $u^A$.
Hence the nodal domains come in pairs and $\mu(u^A)$
 is even. 
\item  
If $u^S$ is a symmetric eigenfunction, then there are two classes of
nodal domains~:
\begin{itemize}
\item  the symmetric domains, 
\begin{equation}\label{DSS}
D_{i,S}=ID_{i,S}\,,
\end{equation}
\item the symmetric pairs of domains $D^1_{i,S}, D^2_{i,S}$ so that 
\begin{equation}\label{DSA}
ID^1_{i,S}=D^2_{i,S}\,. 
\end{equation}
\end{itemize}
\end{enumerate}
\begin{theorem}\label{CourantInv}~\\
Suppose that $\Omega$ satisfies the symmetry assumption \eqref{OVinv}.
Then, if $(\lambda^A_k,u^A)$ is a spectral pair for $H_A(\Omega)$, we have 
\begin{equation}\label{mua}
\mu(u^A)\le 2k\,.
\end{equation}
If $(\lambda^S_k, u^S)$ is a spectral pair for $H_S(\Omega)$ and if we
denote by $\ell(k)$  the number of  pairs of  nodal domains of
$u^S$ satisfying \eqref{DSA} 
and $m(k)$ the number of domains satisfying \eqref{DSS}, then we have  
\begin{equation}\label{muS}
\ell(k)+m(k)\le k\,,
\end{equation}
and  
\begin{equation}\label{muSkl}
\mu(u^S)\le k+\ell(k)\,.
\end{equation}
\end{theorem}
\begin{remark}\label{Icourant}~\\
 Of course the original  Courant Theorem holds, but the above result gives additional informations.
\end{remark}

\textbf{Proof.}~\\
We just have to mimick the proof of Courant's original theorem. Let us
first show \eqref{mua}. We can of course add the condition that~:
$$
\lambda^A_{k-1} < \lambda^A_k\;.
$$
Assume for contradiction that, for some 
 $u^A_k$ we have $\mu(u^A_k)>2k$. To each pair $(D_i,\: ID_i)$ of nodal
domains of $u^A_k$,  we associate the corresponding 
ground states, so that $$H(D_i)\phi_i=\la_k^A\phi_i,\:\phi_i\in
W_0^{1,2}(D_i)\,,$$ and, with $I\phi_i=-\phi_i\circ I\,$,
 $$H(ID_i,V)I\phi_i=\la_k^AI\phi_i\,.$$
We use the variational principle in the form domain of $H_A\,$.
We have 
\begin{equation}\label{lkas}
\la_k^A=\inf_{\varphi^A\bot\mathcal
  Q_A^{k-1}}\frac{\int_\Om\big(|\nabla
  \varphi^A|^2\big)\,d\mu_{\mathbb S^2}}
{\int_\Om|\varphi^A|^2\,d\mu_{\mathbb S^2}}
\end{equation}
where $I\varphi^A=-\varphi^A$ and $\varphi^A\in W_0^{1,2}(\Om)$.  Here
$\mathcal Q_A^{k-1}$
 is just the space spanned by the 
first $(k-1)$  eigenfunctions of $H_A$.  
We  proceed now as  in the proof of Courant's nodal Theorem.  Hence, in
other words, we have just replaced in this proof 
the single domains by pairs of domains.

The proof of \eqref{muSkl} is similar.
\finbox

We can also find some immediate consequences concerning the relation between the $\sigma_A$ and $\sigma_S$. 
Take for instance a spectral pair  $(u^A,\la^A_j)$ and assume that $\mu(u^A)=2k$. Then we can
construct from the $2k$  ground states of each connected component 
$k$ symmetric ones, each one being supported in a symmetric pair of components. By the variational principle, this time for $H_S$ we obtain
\begin{equation}\label{lksym}
\la_k^S=\inf_{\varphi^S\bot\mathcal
  Q_S^{k-1}}\frac{\int_\Om|\nabla\varphi^S|^2\,d\mu_{\mathbb
    S^2}}
{\int_\Om|\varphi^S|^2d\mu_{\mathbb S^2}}.
\end{equation}

This implies~:
\begin{proposition}~\\
 Any eigenvalue $\la^A$ of $H_A$,  whose
corresponding eigenspace contains an eigenfunction with $2k$ nodal 
domains, satisfies $\la^A\ge\la_k^S$. 
\end{proposition} 
A  similar argument can be also made for the symmetric case if $\ell(k)>0$.
This gives us new versions of Courant-sharp properties.\\
If we call {\bf pair symmetric partition}  a partition  which is invariant 
by
the symmetry but such that no element of the
 partition is invariant,  we   have the following Courant-sharp
 analog~:\\
\begin{theorem}\label{CourantA}~\\
If, for some eigenvalue $\lambda_k^A$ of $H_A(\Omega)$, there exists an eigenfunction $u^A$ such
that
 $\mu(u^A) =2 k$, then the corresponding family of nodal domains
 is a minimal pair symmetric partition.
\end{theorem}
Note that, if the labelling  of the eigenvalue (counted as eigenvalue of
$H(\Omega)$)
 is $>2k$),  then it is not a minimal $(2k)$-partition of $\Omega$.
\begin{remark}\label{rembesse}~\\
Let us finally mention as connected result (see for example
\cite{Besse}),
 that if $\Omega$ satisfies \eqref{OVinv}, then $\lambda_2(\Omega)=
 \lambda_1^A(\Omega)$.
\end{remark}

\paragraph{Application}~\\
It is known that the eigenfunctions are the restriction to $\mathbb S^2$ of
the homogeneous harmonic polynomials. Moreover, the eigenvalues are
$\ell (\ell +1)$ ($\ell \geq 0$) with multiplicity $(2\ell +1)$.
 Then,
the  Courant nodal  Theorem says that for a spherical harmonic
$u_\ell$ 
 corresponding to $\ell (\ell +1)$, one should have
\beq
\mu (u_\ell) \leq \ell^2 +1\,.
\eeq
As observed in \cite{Ley1}, one can, using the fact that
\beq
u_\ell (-x) =(-1)^{\ell} u_\ell (x)\,,
\eeq
improve this result by using a variant of  Courant's nodal Theorem with symmetry (see
Theorem \ref{CourantInv}) 
 and this leads to the improvement
\beq\label{impCo}
\mu (u_\ell) \leq \ell (\ell -1) +2\,.
\eeq
Let us briefly sketch the proof of \eqref{impCo}. If $\ell$ is odd,
any eigenfunction is odd with respect to inversion. Hence the
number of nodal domains is even $\mu(
u_\ell) = 2 n_\ell$
 and there are no nodal domains invariant by inversion. Using Courant's
 nodal Theorem  for $H_A$, 
we get with $\ell = 2p+1$ that
$$
\begin{array}{ll}
n_\ell &\leq \sum_{q=0}^{p-1} (2 (2q+1) +1)\;\; +1 \\
& = 2 p (p-1) + 3 p +1\\
& = p (2p +1) +1\\
&= \frac 12 \ell (\ell -1)
 +1\,.
\end{array}
$$
If $\ell$ is even, we can only write
$$
\mu(u_\ell) = 2 n_\ell + p_\ell\,,
$$
where $p_\ell$ is the cardinality of the nodal domains which are
invariant by inversion.\\
 Using Courant's nodal  Theorem for $H_S$, 
we get with $\ell = 2p$, that
$$\begin{array}{ll}
n_\ell + p_\ell &\leq \left(\sum_{q=0}^{p-1} (4q +1)\right) +1\\&
 = 2 p(p-1) + p +1\\& = p (2p -1) +1\\& = \frac 12 \ell (\ell -1)
 +1\,.
\end{array}
$$

Using this improved estimate, we immediately obtain ~:
\begin{proposition}~\\
 The only cases where $u_\ell$ can be Courant-sharp are for $\ell =0$
 and $\ell =1$.\end{proposition}
This proposition has the following consequence~:
\begin{theorem}\label{Unter}~\\
A minimal $k$-partition of $\mathbb S^2$ is nodal if and only if
 $k\leq 2$.
\end{theorem}

Note that in  \cite{Ley1, Ley2} the more sophisticated conjecture 
 (verified for $\ell \leq 6$)  is proposed~:
\begin{conjecture}~\\
$$
\mu (u_\ell) \leq \left\{\begin{array}{lr} \frac 12 (\ell+1)^2& \mbox{
  if }
\ell \mbox{ is odd}\,, \\
 \frac 12 \ell (\ell+2)& \mbox{ if }
\ell \mbox{ is even}\,. \end{array}\right.
$$
\end{conjecture}
\begin{remark}~\\
As indicated by D.~Jakobson to one of us, there is also a
probabilistic version of  this conjecture \cite{NaSo}.
V. N. Karpushkin \cite{Kar} has also  the following bound for the number of
components~:\\
$$
\mu (u_\ell) \leq \left\{\begin{array}{lr} (\ell-1)^2 +2  & \mbox{
  if }
\ell \mbox{ is odd}\,, \\
  (\ell-1)^2 +1 & \mbox{ if }
\ell \mbox{ is even} \,.\end{array}\right.
$$
This is for $\ell$ large slightly better than what we obtained
 with the refined Courant-sharp Theorem. 
Let us also mention the recent paper \cite{ErJaNa} and references therein.
\end{remark}

\begin{remark}\label{remcov}~\\
Considering the Laplacian on the double covering $\mathbb S_{\mathcal C}^2$
of $\ddot {\mathbb S}^2:= \mathbb S^2\setminus \{North,South\}$,
Theorems \ref{CourantInv} and \ref{CourantA}
 hold true where  $I$ is replaced by $\mathcal I$ (introduced in \eqref{defIrev})  corresponding to the map $\phi
\mapsto \phi + 2\pi$. The $\mathcal I$-symmetric eigenfunctions can be identified
 to the eigenfunctions of $H(\mathbb S^2)$ by
 $u_S(x)=u(\pi(x))$  and the restriction
 $H_{\mathcal A}(\mathbb S^2_{\mathcal C})$ of $H(\mathbb
 S^2_{\mathcal C})$ to the $\mathcal
 I$-antisymmetric space leads to a new spectrum, which will be
 analyzed in Section \ref{s8}.
 \end{remark}

\section{On topological properties of minimal $3$-partitions of the
  sphere}
As in the case of planar domains, a classification of the possible
types
 of minimal partitions could simplify the analysis. The case of the
 whole sphere $\mathbb S^2$ shows some difference with for example the
 case of the disk.
\subsection{Around Euler's formula}\label{s3}
As in the case of domains in the plane \cite{HHO},  
 we will use the following  result. 
\begin{proposition}\label{Euler}~\\
Let $\Om$ an open set in $\mathbb S^2$ with piecewise $C^{1,+}$ boundary
 and let  $N\in\mathcal M(\Om)$ such that the associate
 $\mathcal D$  consists of $\mu$ domains $D_1,\dots, D_\mu$. Let $b_0$ be the number of components of 
$\pa \Om$ and $b_1$ be the number of components of $N\cup\pa \Om$. Denote by $\nu(x_i)$ and $\rho(z_i)$
the numbers associated to the $x_i\in X(N)$, respectively  $z_i\in Y(N)$. Then
\begin{equation}\label{Emu}
\mu=b_1-b_0+\sum_{x_i\in X(N)}(\frac{\nu(x_i)}{2}-1)+
\frac{1}{2}\sum_{z_i\in Y(N)}\rho(z_i)+1\,.
\end{equation}
\end{proposition} 
\begin{remark}~\\
In the case when $\Omega =\mathbb S^2$, the statement simply reads
\begin{equation}\label{Emuprime}
\mu=b_1 +\sum_{x_i\in X(N)}(\frac{\nu(x_i)}{2}-1)\;
+1\,,
\end{equation}
where $b_1$ is the number of components of $N$.

\end{remark}
\subsection{Application to $3$- and $4$-partitions.}
\paragraph{The case of $3$-partitions}~\\
Let us analyze in this spirit the topology of minimal $3$-partitions
of $\mathbb S^2$. \\
First we recall that a minimal $3$-partition cannot be nodal. The
multiplicity of the second eigenvalue of $-\Delta_{\mathbb S^2}$ is indeed
$3$.
Hence, our minimal $3$-partition cannot be admissible.\\ Let us look now
to the information given by Euler's formula. We argue like in \cite{HHO}
for the case of the disk.
We recall that at any critical point $x_c$ of $N$, 
\begin{equation}
\nu(x_c)\geq 3\,.
\end{equation}
Hence \eqref{Emuprime} implies that $b_1\leq 2$ and  we have
\begin{equation}
1\leq b_1\leq 2\,.
\end{equation}

When $b_1=1$, we get
as unique solution $\# X(N) =2$. So $X(N)$ consists of two points
$x_1$ and $x_2$ such that  $\nu(x_i)=3$ for $i=1$ and
$2$. The other case when $\# X(N) =1$ leads indeed to an even number
of half lines arriving to the unique critical point and to an
admissible (hence excluded) partition.\\
When  $b_1=2$, we obtain that $X(N)$ is empty and the partition
should be admissible which is excluded. Hence we have shown
\begin{proposition}\label{propotop}~\\
If $\mathcal D$ is a regular non admissible  strong $3$-partition of $\mathbb S^2$, then 
 $X(N)$ consists of two points $x_1$ and $x_2$ such that $\nu(x_i)=3$,
 and 
 $N$  consists of 
 three non crossing (except at their ends) arcs joining  the two points
 $x_1$ and $x_2$.
\end{proposition}
In particular this can be applied to minimal $3$-partitions of
$\mathbb S^2$.
\paragraph{The case of $4$-partitions}~\\
We can analyze in the same way the case of non admissible
 $4$-partitions.  Euler's formula leads to the following classification.
\begin{proposition}\label{propotop4}~\\
If $\mathcal D$ is a regular non admissible  strong $4$-partition of
$\mathbb S^2$, then we are in one of the following cases~:
\begin{itemize}
\item
 $X(N)$ consists of four  points $x_i$ ($i=1,\dots, 4$)  such that $\nu(x_i)=3$,
 and 
 $N$  consists of 
 six non crossing (except at their ends) segments, each one  joining  two points
 $x_i$ and $x_j$ ($i\neq j$).
\item $X(N)$ consists of three points  $x_i$ ($i=1,2, 3$) 
 such that $\nu(x_1)=\nu(x_2)=3$, $\nu(x_3)=4$ and $N$ consists of five
 non crossing (except at their ends)  segments joining two critical points.
\item $X(N)$ consists of two  points $x_i$ ($i=1, 2$)  such that $\nu(x_i)=3$, and 
 $N$  consists of 
 three non crossing (except at their ends) segments joining  the two points
 $x_1$ and $x_2$ and of one closed line.
\item $X(N)$ consists of two  points $x_i$ ($i=1, 2$) such that
 $\nu(x_1)=3$, $\nu(x_2)=5$ and $N$ consists of four 
 non crossing (except at their ends)  segments joining the two critical points
 and of one non crossing (except at his ends) segment starting from one critical point and coming back
 to the same one.
\end{itemize}
\end{proposition}
Note that the spherical tetrahedron corresponds to the first type and
we recall from Theorem \ref{Unter} that minimal $4$-partitions
 are not admissible.
\section{Lyustenik-Shnirelman Theorem and proof of Proposition \ref{antipodal}}\label{s6}

As  we have shown in the previous section $N(\mathcal D)$ consists of two points $x_1\neq x_2$ and $3$ mutually non-crossing arcs $\gamma_1, \gamma_2,\gamma_3$
connecting $x_1$ and $x_2$.  
This means that each $D_i$ has a boundary which is a closed curve
which is away from $x_1,x_2$ smooth.\\

We first recall the well known theorem of Lyusternik and Shnirelman from 1930, that can be found for instance in \cite{M}
on page 23. It states the following. 
\begin{theorem}\label{LS}~\\
 Suppose $S_1,S_2,\dots, S_d$ are closed subsets of $\mathbb S^{d-1}$
 such 
that $\cup_{i=1}^dS_i=\mathbb S^{d-1}$. Then there is at least one $S_i$ that  contains a pair of antipodal points.
\end{theorem}

 We will use this theorem in the case $d=3$  and apply it with  $S_1,S_2,S_3$ defined by 
\begin{equation}\label{SD}
S_i=\overline D_i.  
\end{equation}
In order to prove Proposition \ref{antipodal} it suffices to show that 
\begin{equation}\label{IN}
 N(\mathcal D)\cap  I\, N(\mathcal D)\neq \emptyset\,, 
\end{equation}
where we recall that $I$ is the antipodal map.\\

By Theorem \ref{LS} we know that there is an $S_i=\overline D_i$ which contains a pair of 
antipodal points. After relabelling the $D_i$'s, we can  assume that 
\begin{equation}\label{LS1}
 I \, \overline{ D_1} \cap \overline{ D_1 } \neq \emptyset\;,
\end{equation}
and the goal is to show that
\begin{equation}
 I \, \pa D_1 \cap \pa D_1 \neq \emptyset\;.
\end{equation}
Our $D_i$'s have the properties of $3$-minimal partitions established
in the previous section. In
particular $\pa D_1$ has one component and is also the boundary of
$\pa D_{13}$, where $D_{13}=\Inte ( \overline{D_2} \cup \overline{D_3})$. 

The proof is by contradiction. Let us assume by contradiction that
\begin{equation}\label{contradiction}
 I \, \pa D_1 \cap \pa D_1 = \emptyset\;.
\end{equation}
Then there are two possible cases
\begin{enumerate}
\item[Case a:] $ I\, \pa D_1 \subset D_1$
\item[Case b:]   $ I\,\pa D_1 \cap \overline{ D_1}=\emptyset $
\end{enumerate}

{\bf Let us start with Case a.} Again there are two possibilities.
\begin{enumerate}
\item[Case a1:] $I\, D_1 \subset\subset D_1$ ($\subset\subset$ means
 compactly included).
\item[Case a2:] $I\, D_2 \subset\subset D_1$
\end{enumerate}
But Case $a1$ is in contradiction with the fact that $I$ is an
isometry and 
Case~$a2$ is in contradiction with $\lambda (D_2)=\lambda(D_1)$. 
Of course we could have taken $D_3$ instead of  $D_2$.

{\bf We now look at Case b.} Then $ I \, \pa D_1$ is delimiting in $\mathbb S^2$ two components and
$D_1$ is compactly supported in one component. One of the component is
 $I\, D_1$. But \eqref{LS1} implies that $D_1$ is in this last component~:
$$
D_1 \subset\subset I \, D_1\,.
$$
This can not be true for two isometric domains. Hence we have a
contradiction  with \eqref{contradiction} in all the cases. This
achives the proof of the proposition.

\section{The Laplacian on $\mathbb S^2_{\mathcal C}$} \label{s8}
\subsection{Spherical harmonics with half integers}
These spherical harmonics appear from the beginning
 of Quantum mechanics in connection with  the representation theory \cite{Pau}. We
 refer to \cite{Fl} (Problem 56 (NB2) in the first volume together
 with Problem 133 in the second volume). We are looking for
 eigenfunctions
 of the Friedrichs extension of 
\beq
{\bf L}^2= -\frac{1}{\sin^2\theta} \frac{\pa^2}{\pa \phi^2} - \frac{1}{\sin
 \theta} \frac{\pa}{\pa \theta} \, \sin \theta \, \frac{\pa}{\pa \theta}
\eeq
in $L^2(\sin \theta d\theta \,d\phi)$, 
satisfying
\beq
{\bf L}^2 Y_{\ell m}  = \ell (\ell +1)  Y_{\ell m}\,.
\eeq
The standard spherical harmonics, corresponding to $\ell\geq 0$
 are defined,  for an integer $m \in \{-\ell,\dots,\ell\}$, by
\beq\label{YLM}
Y_{\ell m}(\theta,\phi) = c_{\ell,m} \exp i m \phi \frac{1}{\sin^m\theta}
 (-\frac{1}{\sin \theta} \frac{d}{d\theta})^{\ell -m} \sin^{2\ell} \theta \,,
\eeq
where $c_{\ell,m}$ is an explicit normalization constant.\\
For future extensions, we prefer to take this as a definition for
 $m\geq 0$ and then to observe that
\beq\label{YLMa}
Y_{\ell, -m} = \hat c_{\ell ,m}\overline{ Y_{\ell,m}}\,.
\eeq
For $\ell = 0$, we get $m=0$ and the constant. For $\ell =1$, we
 obtain, for $m=1$, the function  $(\theta,\phi)\mapsto \sin \theta
 \exp i \phi$ and for $m=-1$, the function 
 $\sin \theta \exp -i \phi$ and for $m=0$ the function  $\cos \theta$, which shows
 that the multiplicity is $3$ for the eigenvalue $2$.\\
Of course concerning nodal sets, we look at the real valued functions\break
  $(\theta,\phi)\mapsto \sin \theta
 \cos \phi$ and  $(\theta,\phi)\mapsto \sin \theta
 \sin \phi$ for $|m|=1$.

As  observed a long time ago, these formulas still define
eigenfunctions
 for pairs $(\ell, m)$ with $\ell$ a positive half-integer (and not integer),
  $m \in  \{-\ell,\dots,\ell\}$ and $m-\ell$ integer.\\
For definiteness, we prefer (in the half-integer case) to only consider the pairs with $\ell >0$ and
$m>0$
 and to complete the set of eigenfunctions by introducing
\begin{equation}\label{YLMb}
\widehat Y_{\ell,m} = \overline{Y_{-\ell,m}}
\end{equation}
These functions are only defined on the double covering $\mathbb
S^{2}_\mathcal C$ of\break $\ddot{ \mathbb
S}^2:=\mathbb S^2 \setminus
\{\{\theta=0\}
 \cup \{\theta =\pi\}\}$, which can be  defined by extending $\phi$ to
 the interval $]-2\pi,2\pi]$.\\
When restricted to $\phi \in ]-\pi,\pi]$, they correspond to the
antiperiodic problem with respect to period $2\pi$ in the $\phi$
variable.\\

To show the completeness it is enough to show that, for given $m>0$,
 the orthogonal family (indexed by $\ell \in \{m+\mathbb N\}$) of functions  $\theta\mapsto \psi_{\ell,m}(\theta):=\frac{1}{\sin^m \theta} (-\frac{1}{\sin
 \theta} \frac{d}{d\theta})^{\ell -m} \sin^{2\ell}\theta$
 span all $L^2(]0,\pi[,\sin \theta \,d\theta)$.\\
For this, we consider $\chi \in C_0^\infty(]0,\pi[)$
 and assume that
$$
\int_{0}^\pi \chi (\theta) \psi_{\ell,m}(\theta) \sin \theta d\theta
 =0\,,\, \forall \ell \in\{ m+ \mathbb N\}\;.
$$
We would like to deduce that this implies $\chi=0$. After a change of
variable
 $t=\cos \theta$ and an integration by parts, we obtain that this problem is equivalent to the
problem to show that, if
$$
\int_{-1}^1 \psi (t ) \;[(1-t^2)^\ell]^{(\ell -m)} \; dt
 =0\,,\, \forall \ell\in \{m+\mathbb N\}\;,
$$
then $\psi =0$.\\
Observing that the space spanned  by the functions
 $(1-t^2)^{-m}  ((1-t^2)^\ell)^{(\ell -m)}$ (which are actually
polynomials
 of exact order $\ell$) 
 is the space of all polynomials we can conclude the completeness.\\ 
 Hence we have obtained the
\begin{theorem}~\\
The spectrum of
 the  Laplace Beltrami operator on $\mathbb S^2_{\mathcal C}$ can be
 described by the eigenvalues $\mu_\ell=\ell (\ell+1)$ ($\ell \in \mathbb N/2$),
 each eigenvalue being of multiplicity $(2\ell +1)$. Moreover the
 $Y_{\ell,m}$, as introduced in \eqref{YLM}, \eqref{YLMa} and
 \eqref{YLMb},  define an orthonormal basis of
 the eigenspace $E_{\mu_\ell}$.
\end{theorem}

In particular, for $\ell =\frac 12$, we get a basis of two orthogonal real 
 eigenfunctions $\sin \frac{\phi}{2}
(\sin \theta)^{\frac 12}$  and $\cos \frac \phi 2 (\sin \theta)
 ^{\frac 12}$
 of the eigenspace associated with $\frac 34$. 
 For $\ell =\frac 32$, the multiplicity is  $4$
 and the functions 
 $\sin \frac{3 \phi}{2}(\sin \theta)^{\frac 32}$, 
 $\cos \frac{3 \phi}{2}(\sin \theta)^{\frac 32}$, $\sin \frac \phi 2 
(\sin\theta)^{\frac 12}  \cos \theta$ and  $\cos \frac \phi 2 (\sin
 \theta)^{\frac 12} \cos \theta$ form  a basis of the eigenspace 
 associated with the eigenvalue $\frac{15}{4}$.
\subsection{Covering argument and minimal partition}
Here we give one part of the proof of Proposition~\ref{centera}.
\begin{lemma}\label{lemmacentera}~\\
Let us assume that there exists a $3$-minimal partition $\mathcal D=(D_1,D_2,D_3)$ of $\mathbb S^2$ containing two antipodal
 points in its boundary.  Then,  considering the associated punctured
 $\ddot{\mathbb S}^2$ and the corresponding double covering 
  $\mathbb S^2_{\mathcal C}$ and  denoting  by $\Pi$ the projection
of $\mathbb S^2_{\mathcal C}$ on $\ddot{\mathbb S}^2$, $\Pi^{-1}(D_i)$
 consists of two components and $\pi^{-1} (\mathcal D)$ defines a
 $6$-partition of $\mathbb S^2_{\mathcal C}$ which is pairwise symmetric.
\end{lemma}

The only point to observe is that, according to the property of a
minimal partition established in Proposition \ref{propotop},  the boundary of the partition
necessarily contains a ``broken'' line joining the two antipodal
points.\\
Using the minimax principle, one immediately gets that, {\bf under the
assumption of the lemma},
\begin{equation}\label{inega}
 \mathfrak L_3(\mathbb S^2) \geq \lambda^3_{AS}(\mathbb S^2_{\mathcal
  C})\,,
\end{equation}
where $
\lambda^3_{AS}$ is the third eigenvalue of the Laplace-Beltrami
operator on $\mathbb S^2_{\mathcal
  C}$ restricted to the antisymmetric spectrum.
$\lambda^3_{AS}(\mathbb S^2_{\mathcal
  C})$ will be computed in the next subsection.

\subsection{Covering argument and Courant-sharp eigenvalues}
In the case of the double covering $\mathbb S^2_{\mathcal C}$ of $\ddot{\mathbb S}^2$,
 we have seen that we have to add the antisymmetric (or antiperiodic)
 spectrum  (corresponding to the map $\Pi$, which writes  in spherical
 coordinates~: 
       $\phi \mapsto \phi+2\pi$). This  adds
 the eigenvalue $\frac 34 = \frac 12 (1 +\frac 12)$ with multiplicity
       $2$
 and the eigenvalue $\frac {15}{4} = \frac 32 (1+\frac 32)$ 
with multiplicity $4$.
Hence $\frac {15}{4}$ is the $7$-th eigenvalue of the Laplacian on
 $\mathbb S^2_{\mathcal C}$ 
 hence not Courant-sharp, but it is the third antisymmetric eigenvalue
$$\lambda^3_{AS} = \frac {15}{4}\,.$$
Hence, observing that the nodal set of an eigenfunction associated
 to $\lambda^3_{AS}$ has six nodal domains which are pairwise
 symmetric
 and giving by projection the $Y$-partition, we immediately obtain
 that under the assumption of the lemma
$$
\mathfrak L_3(\mathbb S^2) = \frac{15}{4}\;.
$$

But Proposition \ref{centera} says more. For getting this result, we
have to prove the following proposition~:
\begin{proposition}\label{last}~\\
Let $\mathfrak L_{2\ell}^{AS}(\mathbb S^2_{\mathcal C})$ the infimum
 obtained over the pairwise symmetric (by $\Pi$) $(2 \ell)$- partitions 
 of $\mathbb S^2_{\mathcal C}$. Then, if 
 $$\mathfrak L_{2\ell}^{AS}(\mathbb S^2_{\mathcal C})=
 \lambda^\ell_{AS}\,,
$$
then $ \lambda^\ell_{AS}$ is Courant-sharp in the sense of the
antisymmetric spectrum and any minimal pairwise symmetric 
$(2\ell)-$partition  is nodal.
\end{proposition}
The proof is the same as for Theorem 1.17  in \cite{HHOT} and Theorem
2.6 in \cite{HHO}, with the
difference that we consider everywhere antisymmetric states.\\

Applying this proposition for $\ell =3$, we have the proof of
Proposition \ref{centera}.

\section{On Bishop's approach for   minimal $2$-partitions and
  extensions
 to strong $k$-partitions}\label{s5}
\subsection{Main result for $k=2$}
For $2$-partitions, it is 
 immediate to show that  the minimal  $2$-partitions realizing  
$\mathfrak L_{2}(\mathbb S^2)$  are  given by the two
 hemispheres. One is indeed in the Courant-sharp situation. 
The case of $\mathfrak L_{2,p}(\mathbb S^2)$ for $p<\infty$
 is more difficult. Bishop has described in \cite{Bis} how one can
 show
 that the
 minimal $2$-partitions realizing  $\mathfrak L_{2,1}(\mathbb S^2)$
 are  also 
given by the two
 hemispheres.  It is then easy to see that it implies the property
 for any $p\in [1,+\infty[$. Hence we obtain~:
\begin{theorem}\label{ThBishop}~\\
For any $p\in[1,+\infty]$, $\mathfrak L_{2,p}(\mathbb S^2)$
 is realized by the partition of $\mathbb S^2$ by two hemispheres.
\end{theorem}
The proof is based on two theorems due respectively to Sperner \cite{Sper}
 and Friedland-Hayman \cite{FrHa}. We will discuss their proof because it  will
 have some consequences for the analysis of 
 minimal 
$3$ and $4$-partitions.\\
\subsection{The lower bounds of Sperner and Friedland-Hayman}
For a given domain $D$  on the unit sphere $S^{m-1}$ in $\mathbb R^m$,
 Sperner shows
 the following theorem, which plays on the sphere the same role
as the Faber-Krahn Inequality  plays in $\mathbb R^m$~:
\begin{theorem}\label{Th5.1}~\\
Among all sets $E\subset \mathbb S^{m-1}$ with given $(m-1)$-dimensional surface area
$\sigma_m S$ on  (with $\sigma_m$
being the area
 of  $\mathbb S^{m-1}$), a spherical cap has the smallest characteristic constant.
\end{theorem}
Here the characteristic constant for a domain $D$ is related to the
ground state energy by 
\begin{equation}
\lambda (D) = \alpha (D) (\alpha(D) +m-2)\,,
\end{equation}
with $\alpha (D)\geq 0\,$.\\
We introduce for short 
\begin{equation}
\alpha (S,m)=\alpha (\mathcal{SC}(\sigma_mS))\,,
\end{equation}
where $\mathcal{SC}(\sigma_mS)$ is a spherical cap of surface area
$\sigma_m S$.

This theorem is not sufficient in itself for the problem. The second
 ingredient\footnote{See for example \cite{AltCaf} p.~441 or C.~Bishop \cite{Bis}}
 is a lower bound of $\alpha (S,m)$ by various convex decreasing functions. It is
 indeed proven\footnote{We only write the result for $m=3$ but
   \eqref{prem}
 holds for any $m\geq 3$.}
 in \cite{FrHa} that~:
\begin{theorem}\label{Th4}~\\
We have the  following lower bound ~: 
\begin{equation}\label{prem}
\alpha (S,3)\geq \Phi_3(S)\,,\; 
\end{equation}
where $\Phi_3$ is the convex decreasing function
 defined by 
\beq
 \Phi_3 = \max ( {\widehat \Phi_3}, \Phi_\infty) \,,
\eeq
\begin{equation}
\Phi_\infty (S)=\left\{ 
\begin{array}{ll}
 \frac 12\log (\frac{1}{4S})+ \frac 32\,, &\mbox{ if } 0<S\leq \frac 14\,,\\
2 (1-S)\,,&\mbox{ if } \frac 14 \leq S<1\,.\end{array}\right. 
\end{equation}
\beq 
\widehat \Phi_3(S) = \left\{ 
\begin{array}{ll}
2 (1-S)\,,&\mbox{ if } \frac 12 \leq S <1\,,\\
 \frac 12 j_0\left( \frac 1S -
 \frac 12\right)^{\frac 12}-\frac 12 \,,&\mbox{ if }  S <\frac 12\,,
\end{array}\right. 
\eeq
and  $j_0$ being the first zero of Bessel's function of order $0$~:
\beq 
j_0\sim  2.4048\,.
\eeq
\end{theorem}
\subsection{Bishop's proof for $2$-partitions}
With these two ingredients, we observe (following a remark of
C.~Bishop) that, for a $2$-partition, we have
necessarily
\begin{equation}
\alpha (D_1)+ \alpha (D_2) \geq 2\,,
\end{equation}
the equality being obtained for two hemispheres.\\
The minimization for the sum corresponds to
\begin{equation}
\inf \left(\alpha (D_1) (\alpha(D_1) +1) +  \alpha (D_2) (\alpha(D_2)
+1)\right)\,.
\end{equation}
This infimum is surely larger or equal to
$$
\inf_{\alpha_1+\alpha_2\geq 2\,,\; \alpha_1\geq 0,\alpha_2\geq 0}
\left(\alpha_1(\alpha_1+1) +\alpha_2 (\alpha_2 +1)\right)\,.
$$
It is then easy to see that the infimum is obtained for
 $\alpha_1=\alpha_2 =1$,
\begin{equation}\label{cas2}
\inf_{\alpha_1+\alpha_2\geq 2\,,\; \alpha_1\geq 0,\alpha_2\geq 0}
\left(\alpha_1(\alpha_1+1) +\alpha_2 (\alpha_2 +1)\right)=4\,.
\end{equation}
 This gives a lower bound for $\mathfrak L_{2,1} (\mathbb S^2)$
 which is equal to the upper bound of  $\mathfrak
 L_2(\mathbb S^2)$ and which is attained for the two hemispheres. This
 achieves
 the proof of Theorem \ref{ThBishop}.
\begin{remark}\label{remsusanna}~\\
A natural question is to determine under which condition the infimum
 of $\Lambda^1(\mathcal D)$ for $\mathcal D \in \mathfrak O_2$ 
 is realized for a pair $(D_1,D_2)$ such that $\lambda(D_1)=\lambda(D_2)\,$. Let us illustrate the
 question by a simple example.
 If we consider two disks
 $C_1$ and $C_2$ such that $\lambda(C_1)<\lambda(C_2)\leq
 \lambda_2(C_1)\,$,
 it is not too difficult to see that if we take $\Omega$ as the union
 of these two disks and of a thin channel joining the two disks, then
 $\mathfrak L_2(\Omega)=\lambda_2(\Omega)$ 
 will be very close to $\lambda(C_2)$
 and the infimum of  $\Lambda^1(\mathcal D)$ will be less than $\frac 12
 (\lambda(C_1)+\lambda(C_2))\,$. Hence we will have strict inequality 
 if the channel is small enough. We refer to \cite{BHM1,BHM2,Ann,JM}
 for the spectral analysis of this type of situation. These authors
 are actually more interested in the symmetric situation where
 tunneling plays an important role.
\end{remark}
\subsection{Application to general $k$-partitions}\label{ss53}
One can also discuss what can be obtained in the same spirit 
for  $k$-partitions ($k\geq 3$). This will not lead to the proof
 of Bishop's conjecture but give rather accurate lower bounds
 corresponding in a slightly different context to the ones proposed in
 Friedland-Hayman \cite{FrHa} for harmonic functions in cones of
 $\mathbb R^m$.\\

Let us first mention the easy result extending \eqref{cas2}.
\begin{lemma}\label{lemmeeasy}~\\
Let $k\in \mathbb N^*$ and $\rho >0$. 
If $$
T^{k,\rho}:=\{\alpha \in \overline{\mathbb R}_{+}^{\;k}\;\big |\;\sum_{j=1}^k 
\alpha_j\geq \rho\},
$$
then
$$
\frac 1k \inf_{\alpha\in T^{k,\rho}} \sum_{j=1}^k \alpha_j (\alpha_j +1)
 \geq \frac\rho k (\frac \rho k +1)\,.
$$
\end{lemma}

For a $k$-partition $\mathfrak D=(D_1,\dots,D_k)$, the corresponding characteristic numbers satisfy~:
\beq
\sum_j\alpha (D_j)
 \geq \sum_j \Phi_3 (S_j)
\eeq
with 
\beq
 \sigma_3\, S_j = \mbox{ Area } (D_j)\; (j=1,\,\dots,\,k) \,,
\eeq
and 
\beq
\sum_j S_j =1\,.
\eeq
Using the convexity of $\Phi_3$, we obtain
\begin{proposition}~\\
If $\mathfrak D= (D_j)_{j=1,\dots,k}$ is a strong $k$-partition
 of $\mathbb S^2$, then
\begin{equation}\label{lbalpha}
\frac 1k\, \sum_{j=1}^k \alpha(D_j)\geq  \Phi_3(\frac 1k)\,.
\end{equation}
\end{proposition}

Applying Lemma \ref{lemmeeasy} with  $\rho = k\Phi_3(\frac 1k) $, this leads
 together with \eqref{lbalpha} and \eqref{monoto} to the lower bound
 of $\mathfrak L_{k,1}(\mathbb S^2)$~:
\begin{proposition}\label{prop4.6}
\begin{equation}
\mathfrak L_k (\mathbb S^2) \geq \mathfrak L_{k,1}(\mathbb S^2) 
\geq \Phi_3(\frac 1k)\,( 1 + \Phi_3(\frac 1k)) \,.
\end{equation}
\end{proposition}
Let us see what it gives coming back to the definition of $\Phi_3$.
\begin{corollary}
\begin{equation}
\mathfrak L_k (\mathbb S^2) \geq \mathfrak L_{k,1}(\mathbb S^2) 
\geq \gamma_k  \,,
\end{equation}
with
\beq
\gamma_k:= \Phi_\infty(\frac 1k)\,( 1 + \Phi_\infty(\frac 1k))\,,
\eeq
and
\beq
 \Phi_\infty(\frac 1k):=\left\{ 
\begin{array}{ll}  \frac{2(k-1)}{k} &\mbox{ if } k \leq  4\,,\\
  2 \log (\frac k4)  + \frac 32 & \mbox{ if } k>4\,.
\end{array}\right.
\eeq
In particular
\beq
\gamma_2=2\,,\, \gamma_3 = \frac{28}{9}\,,\, \gamma_4= \frac{15}{4}\;.
\eeq
\end{corollary}
We note that $\gamma_2$ is optimal and  that $\gamma_3 <
\frac{15}{4}$. Hence for $k=3$, the lower bound is not optimal and
does not lead to a proof of Bishop's Conjecture. 
Let us now consider the estimates associated with $\Phi_3$.
\begin{corollary}
\begin{equation}
\mathfrak L_k (\mathbb S^2) \geq \mathfrak L_{k,1}(\mathbb S^2) 
\geq \delta_k  \,,
\end{equation}
with
\beq
\delta_k:= \widehat \Phi_3(\frac 1k)\,( 1 + \widehat \Phi_3(\frac 1k))\,,
\eeq
In particular
\beq
 \delta_3 = \frac 5 8  j_0^2  -\frac 14\,,\, 
\delta_4= \frac 78 j_0^2 -\frac 14 \;.
\eeq
\end{corollary}
\subsection{Discussion for the cases $k=3$ and $k=4$.}
We observe that $\delta_3 >\gamma_3$ and $\delta_4 >\gamma_4$.
Small computations\footnote{already done in \cite{FrHa}} show indeed  that 
\begin{equation}\label{valimpa}
\widehat \Phi_3(\frac 13) \sim 1.401\,,
\end{equation}
which is higher than $\Phi_\infty (\frac 13)=\frac 43$, 
and
\begin{equation}
\widehat \Phi_3(\frac 14)\sim 1.748\,,
\end{equation}
which is higher than $\Phi_\infty(\frac 14)= \frac 32\,$.
This 
leads to  the lower bound~:
\begin{proposition}
\begin{equation}\label{minoL4b}
\mathfrak L_4 (\mathbb S^2) \geq \mathfrak L_{4,1}(\mathbb S^2) >
15/4 =\mathfrak L_3 (\mathbb S^2)\,.
\end{equation}
\end{proposition}
In particular the best lower bound of $\mathfrak L_{4,1}(\mathbb S^2) $
 is approximately
\begin{equation}\label{minoL4app}
\delta_4 \sim 4.8035\,.
\end{equation}
Note that by a third method, one can find in \cite{FrHa} (Theorem 5 and
table~1, p.~155, computed by J.G. Wendel) another convex function
$\tilde \Phi$ such that $\alpha(D)\geq \tilde \Phi(S)$ and
\begin{equation}\label{valimpabis}
\tilde \Phi (\frac 13)\sim 1.41167\,.
\end{equation}
Note that $\tilde \Phi(\frac 14) < \widehat \Phi_3(\frac14)$ so this improvment
occurs only for $3$-partitions.\\

In the case of $\mathbb S^2$, unlike the case of the square or the
disk, the minimal $4$-partition is not nodal (as proven in Theorem \ref{Unter}).  Note that this implies that the $4$-minimal partition
realizing $\mathfrak L_{4,p}(\Omega)$
 for  $p\in [1,+\infty]$ is
 neither nodal.
\\
 As already mentioned in \cite{FrHa}, there is at least a natural candidate which is the
spherical regular tetrahedron. Numerical
computations\footnote{transmitted to us by M. Costabel}  give, for
the corresponding $4$-partition $\mathcal D_4^{Tetra}$, 
\begin{equation}
\Lambda (\mathcal D_4^{Tetra}) \sim 5.13\,.
\end{equation}
Hence we obtain   that
\begin{equation}\label{Tetra}
\frac{15}{4} < \mathfrak L_4(\mathbb S^2) \leq \Lambda (\mathcal D_4^{Tetra})
< 6 =L_4(\mathbb S^2)\,.
\end{equation}
It is interesting to compare it with \eqref{minoL4app}.\\
According to a personal communication of M. Dauge, one can also
observe
 that the largest circle inside a face of the tetrahedron is actually
 a nodal line corresponding to an eigenfunction with eigenvalue
 $6$ (up to a rotation, this is the (restriction to $\mathbb S^2$ of the)
 harmonic polynomial  $\mathbb S^2\ni (x,y,z)\mapsto x^2+y^2 - 2 z^2\,$. This
 gives  directly
 the comparison $\Lambda (\mathcal D_4^{Tetra})
< 6\,$.

\subsection{Large $k$ lower bounds}

We can push the argument by looking at the asymptotic as $k\ar +\infty$
of $\delta_k$.
This gives~:
\beq\label{lbk}
\mathfrak L_{1,k}(\mathbb S^2) \geq \frac 14 j_0^2 k  - \frac 18\,
j_0^2 
-\frac 14 \;.
\eeq
We note that at least for $k$ large, this is much better than
 the trivial lower bound
\beq
\mathfrak L_{1,k}(\mathbb S^2) \geq
 \frac{1}{k} \mathfrak L_{k}(\mathbb S^2)\;.
\eeq
We discuss below in Remark \ref{remfermio} an independent improvment.\\
Multiplying \eqref{lbk} by $4\pi$, the area of $\mathbb S^2$
and dividing by $k$, we obtain
\beq
\mbox{Area}(\mathbb S^2)\;\liminf_{k\ar +\infty} \frac{\mathfrak L_{1,k}(\mathbb S^2)}{k}
 \geq \pi j_0^2\,.
\eeq
But $\pi j_0^2$ is the groundstate energy $\lambda(D^1)$
of the Laplacian on the  disk $D^1$ in $\mathbb R^2$
 of area $1$. Although not
written explicitly\footnote{The authors mention only the lower bound
  for $\mathfrak L_k(\Omega)$} in \cite{HHOT, BHV:2007}, the Faber-Krahn Inequality gives
for planar
 domains 
\beq
|\Omega| \frac{\mathfrak L_{k,1}(\Omega)}{k} \geq \lambda(D^1)=\pi j_0^2\;.
\eeq
We have not verified the details, but we think that as in
\cite{BHV:2007} 
 for planar domains, we will have
\beq
\mbox{Area}(\mathbb S^2)\;\limsup_{k\ar +\infty} \frac{\mathfrak
  L_{k}(\mathbb S^2)}{k}\leq \lambda(Hexa^1)\,,
\eeq
where $Hexa^1$ denotes the regular hexagon of area $1$.\\

As for the case of plane domains, it is natural to conjecture
(see for example \cite{BHV:2007, CL1} but we first heard of this
question from 
M. Van den Berg five years ago
)
 that~:
\begin{conjecture}~\\
$$
\lim_{k\ar + \infty}  \frac{\mathfrak
  L_{k}(\mathbb S^2)}{k} = \lim_{k\ar + \infty}  \frac{\mathfrak
  L_{k,1}(\mathbb S^2)}{k}=\lambda(Hexa^1)\;.
$$
\end{conjecture}
The first equality in the conjecture corresponds to the idea, which is
well illustrated in the recent paper by Bourdin-Bucur-Oudet \cite{BBO}
 that, asymptotically as $k\ar +\infty$, a minimal $k$-partition for $\Lambda^p$ will
 correspond  to $D_j$'s such that the $\lambda(D_j)$ are equal.

\begin{remark}\label{remfermio}~\\
If $\Omega$ is a regular bounded open set in $\mathbb R^2$ or 
 in $\mathbb S^2$, then 
\begin{equation}\label{fermionic}
\frac 1 k \sum_{j=1}^k \lambda_j(\Omega) \leq \mathfrak L_{k,1}(\Omega)\,.
\end{equation}
The proof is "fermionic". It is enough to apply the minimax characterization
 for the groundstate energy $\lambda^{Fermi,k}$
 of the Dirichlet realization of the  Laplacian on $\Omega^k$ (in $(\mathbb R^2)^k$ or in $(\mathbb S^2)^k$)
 restricted to the Fermionic space $\wedge^k L^2(\Omega)$
 which is $$
 \lambda^{Fermi,k}=\sum_{j=1}^k \lambda_j(\Omega)\,.
$$ For any $k$-partition $\mathcal D$ 
 of $\Omega$, we can
 indeed
 consider the Slater determinant of the normalized groundstates $\phi_j$  of each
$D_j$ and observe that the corresponding energy is
 $k\Lambda^1(\mathcal D)$.\\
This suggests the following conjecture (which is proven for $p=1$
 and $p=+\infty$)~:
\begin{equation}\label{fermionicp}
\left(\frac 1 k \sum_{j=1}^k \lambda_j(\Omega)^p\right)^{\frac 1 p}
 \leq \mathfrak L_{k,p}(\Omega)\,,\,\forall k\geq 1\,,\, \forall p \in [1,+\infty]\,.
\end{equation}
The case when $\Omega$ is the union of  two disks considered in Remark \ref{remsusanna} gives an example
 for $k=2$ where \eqref{fermionic} becomes an equality. In this case,
 we have indeed $\lambda_1(\Omega)=\lambda(C_1)$ and $\lambda_2(\Omega)=\lambda(C_2)$.
  \end{remark}

\newpage 

{\bf Acknowledgements.}~\\
This work was motivated by questions of  A.~Lemenant
 about Bishop's Conjecture.
Discussions with (and numerical computations of) M.~Costabel and
 M.~Dauge were also quite helpful. Many thanks also to A.~Ancona and 
  D.~Jakobson for
 indicating to us useful references.

\footnotesize
\bibliographystyle{plain}

\begin{thebibliography}{1}

\bibitem{AS} M. Abramowitz and I. A. Stegun.
\newblock {\em Handbook of mathematical functions}, 
\newblock Volume
 55 of Applied Math Series. National Bureau of Standards, (1964).


\bibitem{AltCaf} H.W. Alt, L.A. Caffarelli and A. Friedman.
\newblock Variational problems with two phases and their free
boundaries.
\newblock TAMS 282 (2) (1984), p.~431-461.

\bibitem{Ann} C. Ann\'e.
\newblock A note on the generalized Dumbbell Problem.
\newblock Proc. of the AMS 123 (8) (1995), p.~2595-2599.

\bibitem{Bers:1955}
L.~Bers.
\newblock Local behaviour of solutions of general linear equations.
\newblock Commun. Pure Appl. Math. 8 (1955), p.~473-496.

\bibitem{Besse} A. Besse.
\newblock{\em  Manifolds all of whose geodesics are closed.} With appendices
by D.B.A. Epstein, J.-P. Bourguignon, L. B\'erard-Bergery, M. Berger
and J. L. Kazdan.
\newblock  Ergebnisse der Mathematik und ihrer Grenzgebiete,
 93. Springer-Verlag, Berlin-New York, 1978. 


\bibitem{Bis} C.J. Bishop.
\newblock Some questions concerning harmonic measure.
\newblock Dahlberg, B. (ed.) et al., Partial Differential equations
 with minimal smoothness and applications.
\newblock IMA Vol. Math. Appl. 42 (1992), p.~89-97.

\bibitem{BHHO} V.~Bonnaillie-No\"el, B. Helffer and
  T. Hoffmann-Ostenhof.
\newblock  
 Aharonov-Bohm Hamiltonians, isospectrality and  minimal partitions.
\newblock Preprint 2008. To appear in Journal of Physics A.

\bibitem{BHV:2007} V.~Bonnaillie-No\"el, B. Helffer and G. Vial.
\newblock Numerical simulations for nodal domains
 and spectral minimal partitions.
\newblock Preprint 2007. To appear in COCV.

\bibitem{BBO} B.~Bourdin, D.~Bucur, and E.~Oudet.
\newblock Optimal partitions for eigenvalues.
\newblock Preprint January 2009.

\bibitem{BHM1} R.M. Brown, P.D. Hislop, and A. Martinez.
\newblock Lower bounds on the interaction between cavities connected by a thin
tube.
\newblock Duke Math. Journal Vol. 73, No 1 (1994), P.~163-176.

\bibitem{BHM2} R.M. Brown, P.D. Hislop, and A. Martinez.
\newblock Lower bounds on eigenfunctions and the first eigenvalue gap.
\newblock Differential equations with applications to mathematical Physics,
33-49, Math. Sci. Engrg., 192, Academic Press.


\bibitem{CL1} L.A. Caffarelli and Fang Hua Lin.
\newblock An optimal partition problem for eigenvalues.
\newblock Journal of Scientific Computing 31 (1/2) (2006), p. 5-18.

\bibitem{CL2} L.A. Caffarelli and Fang Hua Lin.
\newblock Singularly perturbed elliptic systems and multi-valued
harmonic functions with free boundaries.
\newblock Journal of the AMS 21 (3) (2008), p.~847-862.

\bibitem{CTV0} M. Conti, S. Terracini, and G. Verzini.
\newblock An optimal partition problem related to nonlinear eigenvalues.
\newblock JFA 198  (2003), p.~160-196.

\bibitem{CTV2} M. Conti, S. Terracini, and G. Verzini.
\newblock  A variational problem for the spatial segregation of reaction-diffusion systems.
\newblock Indiana Univ. Math. J. 54 (3) (2005), p. 779-815.

\bibitem{CTV:2005} 
M. Conti, S. Terracini, and  G. Verzini.
\newblock On a class of optimal partition problems related to the Fucik spectrum 
and to the monotonicity formula.
\newblock  Calc. Var. 22 (2005), p.~45-72.
 

\bibitem{ErJaNa} A. Eremenko, D. Jakobson, and N. Nadirashvili.
\newblock On nodal sets and nodal domains on $\mathbb S^2$.
\newblock Annales Institut Fourier 57 (7) (2007), p.~2345-2360.

\bibitem{Fl} S. Fl\"ugge.
\newblock Practical quantum mechanics I, II.
\newblock Die Grundlehren der mathematischen Wissenschaften 177-178.
\newblock Springer Verlag 1971.

\bibitem{FrHa} S. Friedland and W.K. Hayman.
\newblock Eigenvalue Inequalities for the Dirichlet problem on spheres
 and the growth of subharmonic functions.
\newblock Comment. Math. Helvetici 51 (1976), p. 133-161.


\bibitem{He} B. Helffer.
\newblock Domaines nodaux et partitions spectrales minimales.
\newblock S\'eminaire Equations aux D\'eriv\'ees Partielles 2006-2007,
Expos\'e VIII. Publications de l'Ecole Polytechnique.


\bibitem{HH:2005a}
 B.~Helffer, T.~Hoffmann-Ostenhof.
\newblock Converse spectral problems for nodal domains.
\newblock  Moscow Mathematical Journal 7 (2007), p. 67-84.


\bibitem{HHO} B. Helffer, T. Hoffmann-Ostenhof.
\newblock On spectral minimal partitions : new properties 
and applications to the disk.
\newblock Submitted (2008).

\bibitem{HHOT} B. Helffer, T. Hoffmann-Ostenhof, S. Terracini.
\newblock Nodal domains and spectral minimal partitions.
\newblock  Ann. Inst. H. Poincar\'e Anal. Non
  Lin\'eaire  26 (2009), p.~101-138.


\bibitem{Henrot} A. Henrot.
\newblock {\em Extremum Problems for Eigenvalues of Elliptic Operators.}
\newblock Birkh\"auser, Bases, Boston, London, 2006. 


\bibitem{HoMiNa:1999}
T. Hoffmann-Ostenhof, P. Michor, N. Nadirashvili.
\newblock Bounds on the multiplicity of eigenvalues for fixed membranes.
\newblock  GAFA 9 (1999), p. 1169-1188.

\bibitem{JM} S. Jimbo and Y. Morita.
\newblock Remarks on the behavior of certain eigenvalues
 on a perturbed domain with several thin channels.
\newblock CPDE 17 (1992)
, p. 523-552.

\bibitem{Kar} V.N. Karpushkin.
\newblock On the number of components of the complement to some
algebraic curves.
\newblock Russian Math. Surveys 57 (2002), p. 1228-1229.


\bibitem{Ley1} J. Leydold.
\newblock Nodal properties of spherical harmonics.
\newblock PHD 1993 (Vienna University).

\bibitem{Ley2} J. Leydold.
\newblock On the number of nodal domains of spherical harmonics.
\newblock Topology 35 (1996), p. 301-321.


\bibitem{M}
J. Matousek.
\newblock {\em Using the Borsuk-Ulam Theorem. Lectures on Topological 
Methods in Combinatorics and Geometry.}
\newblock Springer 2003.


\bibitem{NaSo} F. Nazarov and M. Sodin.
\newblock On the number of nodal domains of random spherical
harmonics.
\newblock Preprint 2008.

\bibitem{Pau} W. Pauli.
\newblock Helv. Phys. Acta 12  (1939), p.~147.


\bibitem{Pleijel:1956} A.~Pleijel. 
\newblock Remarks on Courant's nodal theorem.
\newblock Comm. Pure. Appl. Math. 9 (1956), p.~543-550. 

\bibitem{Sper} E. Sperner.
\newblock Zur symmetrisierung von Functionen auf Sph\"aren.
\newblock Math. Z. 134 (1973), p. 317-327.



\end{thebibliography}

\scshape 
B. Helffer: D\'epartement de Math\'ematiques, Bat. 425,
Universit\'e Paris-Sud, 91 405 Orsay Cedex, France.

email: Bernard.Helffer@math.u-psud.fr

\scshape
T. Hoffmann-Ostenhof: Institut f\"ur Theoretische Chemie, Universit\"at
Wien, W\"ahringer Strasse 17, A-1090 Wien, Austria and International Erwin
Schr\"odinger Institute for Mathematical Physics, Boltzmanngasse 9, A-1090
Wien, Austria.

email: thoffman@esi.ac.at

\scshape
S. Terracini: 
Universit\`a di Milano Bicocca, Via Cozzi, 53 20125 Milano (Italy).

email: susanna.terracini@unimib.it

\end{document}